\documentclass[12pt,a4paper]{article}
\usepackage{amssymb,graphicx}
\title{A Characteristic Property of \\ Elliptic Pl\"ucker Transformations}
\author{Hans Havlicek}
\date{}

\newtheorem{prop}{PROPOSITION}
\newtheorem{theo}{THEOREM}

\newcounter{zahl}
\newenvironment{axiom}{\begin{list}{\rm \arabic{zahl}.}{\usecounter{zahl}
 \itemsep0em plus 0.1em \parsep0em plus 0.1em \rm}}{\end{list}}

\newcommand{\Ecal}{{\cal E}}
\newcommand{\Fcal}{{\cal F}}
\newcommand{\Lcal}{{\cal L}}

\newcommand{\Pcal}{{\cal P}}
\newcommand{\Rcal}{{\cal R}}
\newcommand{\Scal}{{\cal S}}
\newcommand{\Tcal}{{\cal T}}

\newcommand{\Wcal}{{\cal W}}

\newcommand{\Phat}{\widehat{\cal P}}
\newcommand{\Lhat}{\widehat{\cal L}}

\newcommand{\li}{\parallel_L}
\newcommand{\re}{\parallel_R}
\newcommand{\cp}{\parallel}
\newcommand{\nli}{\,{\not{}\hspace{-0.42em}\parallel}_L\,}
\newcommand{\nre}{\,{\not{}\hspace{-0.42em}\parallel}_R\,}

\newcommand{\Eli}{\Ecal_L}
\newcommand{\Ere}{\Ecal_R}
   \newcommand{\kappali}{{\kappa_L}}
   \newcommand{\kappare}{{\kappa_R}}
   \newcommand{\varphili}{{\varphi_L}}
   \newcommand{\varphire}{{\varphi_R}}
   \newcommand{\varphiliinv}{{\varphi_L^{-1}}}
   \newcommand{\varphireinv}{{\varphi_R^{-1}}}
\newcommand{\Sli}{\Scal_L}
\newcommand{\Sre}{\Scal_R}
\newcommand{\Rli}[2]{\Rcal_L(#1|#2)}
\newcommand{\Rre}[2]{\Rcal_R(#1|#2)}

\newcommand{\rel}{\sim}
\newcommand{\rrel}{\approx}
\newcommand{\nrel}{\not\sim}
\newcommand{\nrrel}{\not\approx}

\newcommand{\PG}[2]{\mbox{$\mbox{{\rm PG}}(#1,#2)$}}

\newcommand{\qu}[1]{\overline{#1}}
\newcommand{\im}{\mbox{\rm im\,}}
\newcommand{\spn}{\mbox{{\rm span\,}}}
\newcommand{\Char}{\mbox{{\rm Char\,}}}
\newcommand\inv{^{-1}}

\newcommand{\abb}[3]{\mbox{$#1\,:\,#2\rightarrow#3$}}
\newcommand{\Abb}[5]{\mbox{$#1\,:\,#2\rightarrow#3,\;#4\mapsto #5$}}
\newcommand{\proof}{{\it Proof. }}
\newcommand{\qed}{\hfill$\Box$}
\newcommand{\zitat}[4]{\bibitem{#1}{\sc #2}: {\sl #3\/}. #4.\vspace{-0.2em}}

\newcommand{\abstand}{\vspace{1em plus0.3em minus 0.3em}}
\begin{document}
\maketitle

%
%
%

\vspace{0.6cm}
{\sl Dedicated to Walter Benz on the occasion of his 65th birthday}

%
%
%
%
\vspace{1.5cm plus0.2cm minus0.2cm}

{\small
We discuss elliptic Pl\"ucker transformations of three-dimensional elliptic
spaces. These are permutations on the set of lines such that any two related
(orthogonally intersecting or identical) lines go over to related lines in
both directions. It will be shown that for ``classical'' elliptic $3$-spaces a
bijection of its lines is already a Pl\"ucker transformation, if related lines
go over to related lines. Moreover, if the ground field admits only surjective
monomorphisms, then ``bijection'' can be replaced by ``injection''.
} 
%
\section{Introduction and Main Results}\label{EIN}

Let $(\Pcal,\Lcal,\pi)$ be a $3$-dimensional {\em elliptic space}, i.e.\ a
projective space $(\Pcal,\Lcal)=\PG3F$ endowed with an elliptic absolute
polarity $\pi$. Points $X,Y\in\Pcal$ are $\pi$-conjugate (orthogonal) if $X\in
Y^\pi$ or, equivalently, if $Y\in X^\pi$. Given a subspace $\Tcal\subset
\Pcal$ we denote by $\Tcal^\pi$ its $\pi$-polar subspace, i.e.\ the set of
points $Y\in \Pcal$ which are $\pi$-conjugate to all points of $\Tcal$. Due to
the absence of $\pi$-self-conjugate ($\pi$-absolute) points, $\Tcal$ and
$\Tcal^\pi$ are complementary subspaces. Moreover, the field $F$ is
necessarily infinite \cite[5.3]{Hi79}.
Let us recall the terminology and some results of \cite{Ha94}: Given $a,b\in
\Lcal$ we put
   \begin{displaymath}\begin{array}{rcll}
   a\rrel b
      & :\Longleftrightarrow & a\cap b^\pi\neq\emptyset \mbox{ and } a\cap
      b\neq\emptyset
      & \mbox{(orthogonally intersecting lines),} \\
   a\rel b
      & :\Longleftrightarrow & a\rrel b \mbox{ or } a=b
      & \mbox{(related lines),} \\
   a\cp b
      & :\Longleftrightarrow & \#\{x\mid x\rrel a\mbox{ and } x\rrel b\}\geq 3
      & \mbox{(Clifford parallel lines).}
   \end{array}\end{displaymath}
The relation $\rrel$ is symmetric and the pair $(\Lcal,\rel)$ is a Pl\"ucker
space in the sense of W.\ Benz \cite{Be92}. An {\em elliptic Pl\"ucker
transformation} is a bijection $\abb{\varphi}{\Lcal}{\Lcal}$ such that
   \begin{equation}\label{PL}
      a\rel b \Longleftrightarrow a^\varphi\rel b^\varphi.
   \end{equation}
In the sequel we shall assume that $(\Pcal,\Lcal,\pi)$ is {\em classical},
i.e., the following conditions hold true:
   \begin{axiom}
      \item The underlying field $F$ is commutative and $\Char F\neq 2$.
      \item $\pi$ is a projective polarity.
      \item There exist Clifford parallel lines $a,b\in\Lcal$ with $b\notin
      \{a,a^\pi\}$.
   \end{axiom}
For example, the real elliptic 3-space fits into this concept. We shall
comment on the third condition in Section \ref{CLIFF}.

In \cite{Ha94} all Pl\"ucker transformations of a classical elliptic $3$-space
have been described in the realm of the ambient space of the Klein quadric
representing the lines of $(\Pcal,\Lcal)$.

The main results of this paper are:
   \begin{theo}\label{TH1}
   Let $(\Pcal,\Lcal,\pi)$ be a $3$-dimensional classical elliptic space.
   If $\abb{\varphi}{\Lcal}{\Lcal}$ is a bijection satisfying
   \begin{equation}\label{TH11}
      a\rel b \Longrightarrow a^\varphi\rel b^\varphi,
   \end{equation}
   then $\varphi$ is an elliptic Pl\"ucker transformation.
\end{theo}

\begin{theo}\label{TH2}
   Let $(\Pcal,\Lcal,\pi)$ be a $3$-dimensional classical elliptic space with
   underlying field $F$. Suppose that there are only surjective monomorphisms
   $F\rightarrow F$. If $\abb{\varphi}{\Lcal}{\Lcal}$ is an injection
   satisfying (\ref{TH11}),
   then $\varphi$ is an elliptic Pl\"ucker transformation.
\end{theo}
We shall establish a series of Propositions in Section \ref{BEW} that end up
in proofs for Theorem \ref{TH1} and Theorem \ref{TH2}.

Let us remark that \cite{Ha94} contains results on generalized elliptic spaces
of dimensions $2$ and $\geq 4$, namely a description of their Pl\"ucker
transformations and characterization in the spirit of Theorem \ref{TH1}. The
corresponding proofs are short and straightforward, whereas the
$3$-dimensional case seems to be much more involved.

For results and references on other groups of Pl\"ucker transformations see,
among others, \cite{Be92}, \cite{Be94}, \cite{Ha94}, \cite{Ha95} and
\cite{Le95}. Finally, we refer to \cite{Bu93}, \cite{FF95}, \cite[p.\
75]{KaKr88}, \cite{KPS93}, \cite{Le54a}, \cite{Le54b}, \cite{Le57} and
\cite{Schr95} for an axiomatic descriptions of polarities, elliptic spaces and
Clifford parallelism as well as a connection with quaternion skew fields.


\section{Clifford Parallelism}\label{CLIFF}

In \cite{Ha94} we have aimed at understanding the line geometry of a classical
elliptic 3-space $(\Pcal,\Lcal,\pi)$ via the ambient space
$(\Phat,\Lhat)=\PG5F$ of the Klein quadric. Write
   \begin{displaymath}
   \Abb{\gamma}{\Lcal}{\Phat}{a}{a^\gamma}
   \end{displaymath}
for the Klein mapping and put $\Gamma:=\Lcal^\gamma=\im\gamma$ for the Klein
quadric. The projective polarity associated with the Klein quadric is named
$\kappa$.

The absolute polarity $\pi$ gives rise to a projective collineation
$\abb{\alpha}{\Phat}{\Phat}$ characterized by $a^{\gamma\alpha} =
a^{\pi\gamma}$ for all $a\in \Lcal$. Since $\pi$ has no self-polar lines,
there is no $\alpha$-invariant point on the Klein quadric. However, since
$(\Pcal,\Lcal,\pi)$ is classical, all $\alpha$-invariant points%
   \footnote{The existence of an $\alpha$-invariant point is equivalent to the
   existence of a $\pi$-invariant general linear complex of lines
   $\subset\Lcal$ or, in other words, is equivalent to the existence of a
   symplectic polarity of $(\Pcal,\Lcal)$ commuting with $\pi$.
   If there would be no $\alpha$-invariant points, then the lines $a$ and
   $a^\pi$ would be the only Clifford parallel lines for any $a\in \Lcal$
   \cite[Lemma 2]{Ha94}.}
form two skew planes of $(\Phat,\Lhat)$, say $\Eli$ and $\Ere$, with
$\Eli^\kappa = \Ere$ \cite[p.\ 45]{Ha94}. We remark that in an appropriate
quadratic extension of \PG3F the absolute polarity $\pi$ becomes the polarity
of a ruled quadric%
   \footnote{If we are given an elliptic projective polarity of \PG3F ($F$
   commutative), then it is possible to obtain self-conjugate points by a
   single quadratic extension of \PG3F, but in general it takes two
   consecutive quadratic extensions to get self-polar lines.}%
.
The two reguli on this quadric go over to distinct conics on the Klein quadric
spanning the planes $\Eli$ and $\Ere$, respectively. Cf.\ part IV of the
fundamental paper by G.\ Wei{\ss} \cite{We78} on real metric line geometry.

We infer from $\Eli\cap \Gamma=\emptyset$ that the polarity of the Klein
quadric induces an elliptic projective polarity in $\Eli$, say $\kappali$,
thus turning $\Eli$ into an elliptic plane. By symmetry of $L={}$``left'' and
$R={}$``right'', this carries over to $\Ere$. The planes $\Eli$ and $\Ere$
give rise to projections
   \begin{eqnarray*}
   {}& \Abb{\lambda}{\Phat\setminus\Eli}{\Ere}{X}{(X\vee\Eli)\cap\Ere}, &{}\\
   {}& \Abb{\rho}   {\Phat\setminus\Ere}{\Eli}{X}{(X\vee\Ere)\cap\Eli}, &{}
   \end{eqnarray*}
with the property that
$(a^{\gamma\lambda}, a^{\gamma\rho}, a^\gamma, a^{\pi\gamma})$
is a harmonic range of points for each $a\in\Lcal$. Let $a,b\in\Lcal$. We
define
   \begin{displaymath}\begin{array}{rcll}
   a\li b
      & :\Longleftrightarrow & a^{\gamma\lambda} = b^{\gamma\lambda}
      & \mbox{(left parallel lines),} \\
   a\re b
      & :\Longleftrightarrow & a^{\gamma\rho} = b^{\gamma\rho}
      & \mbox{(right parallel lines).}
   \end{array}\end{displaymath}
Moreover, by \cite[pp.\ 44--46]{Ha94},
   \begin{eqnarray}
   \label{CP}
   {}&
   a\cp b\Longleftrightarrow  a\li b \mbox{ or }  a\re b,
   &{}\\
   \label{LR}
   {}&
   b\in \{a,a^\pi\} \Longleftrightarrow  a\li b \mbox{ and } a\re b,
   &{}\\
   \label{RREL}
   {}&
   a\rrel b
   \Longleftrightarrow
   a^{\gamma\lambda},b^{\gamma\lambda} \mbox{ $\kappare$-conjugate and }
   a^{\gamma\rho}, b^{\gamma\rho} \mbox{ $\kappali$-conjugate}.
   & {}
   \end{eqnarray}
Left parallelism $\li$ is an equivalence relation. The equivalence class of
$a\in \Lcal$ is an elliptic linear congruence of lines (regular spread)
\cite[Lemma 3]{Ha94}; this spread is denoted by
   \begin{displaymath}
   \Sli(a):=\{x \mid x\li a \}.
   \end{displaymath}
Given a line $p\in\Lcal\setminus \Sli(a)$ then, by the regularity of the
spread $\Sli(a)$,
\begin{displaymath}
   \Rli{a}{p}:=\{x\mid x\li a \mbox{ and } x\cap p\neq\emptyset\}
   \end{displaymath}
is a regulus. These results carry over to $\re$ in an obvious way.

Let $a\rrel b$. If $Q\in\Pcal$, then there exist lines
$x_Q\in\Sli(a)$ and $y_Q\in\Sli(b)$, concurrent at $Q$, since
$\Sli(a)$ and $\Sli(b)$ are spreads. According to \cite[Lemma 5, I]{Ha94}, we
obtain $x_Q\rrel y_Q$, whence
   \begin{equation}
   \label{LILI-SCHNITT}
   \Sli(a)=\{x\in\Sli(a)\mid\exists\, y\in\Sli(b)\mbox{ with }x\rrel y \};
   \end{equation}
$\Sli(b)$ can be described likewise. On the other hand, by \cite[Lemma 5,
II]{Ha94},
   \begin{eqnarray}
   \label{LIRE-SCHNITT}
   {} & \Rli{a}{b} =
   \{x\in\Sli(a)\mid\exists\, y\in\Sre(b)\mbox{ with }x\rrel y \}, & {}\\
   \label{RELI-SCHNITT}
   {} & \Rre{b}{a} =
   \{y\in\Sre(b)\mid\exists\, x\in\Sli(a)\mbox{ with }y\rrel x \}. & {}
   \end{eqnarray}
Consequently, $\Rli{a}{b}$ and $\Rre{b}{a}$ are mutually
opposite reguli%
   \footnote{In a real elliptic $3$-space these two reguli are on the
   well-known {\em Clifford surface}. Formula (\ref{LI-RE}) reflects the fact
   that this quadric admits locally Cartesian coordinates.}
and
   \begin{equation}\label{LI-RE}
   x\rrel y \mbox{ for all } x\in\Rli{a}{b} \mbox{ and all } y\in\Rre{b}{a}.
   \end{equation}
These results remain true if the terms ``left'' and ``right'' are
interchanged.


\section{Proofs}\label{BEW}

In the subsequent Propositions let $(\Pcal,\Lcal,\pi)$ be a $3$-dimensional
classical elliptic space. Suppose, furthermore, that
$\abb{\varphi}{\Lcal}{\Lcal}$ is an injection satisfying (\ref{TH11}).
\begin{prop}\label{A}
   For all $a,b\in \Lcal$ the following properties hold true:
   \begin{eqnarray}
   \label{A0}
   {} & a\rrel b \Longrightarrow a^\varphi\rrel b^\varphi, &{}\\
   \label{A1}
   {} & a\cp b \Longrightarrow a^\varphi\cp b^\varphi,  &{}\\
   \label{A2}
   {} & a^{\pi\varphi} = a^{\varphi\pi},  &{}\\
   \label{A3}
   {} & a\cp b \mbox{ and } b\notin\{a,a^\pi\} \Longrightarrow
   \mbox{ either } a^\varphi\li b^\varphi \mbox{ or } a^\varphi\re b^\varphi.
    &{}
   \end{eqnarray}
\end{prop}
\proof
We deduce (\ref{A0}) from the injectivity of $\varphi$ and (\ref{TH11}). Now
(\ref{A1}) is immediate from the definition of $\cp$, the injectivity of
$\varphi$ and (\ref{A0}). In order to establish (\ref{A2}) choose points
$A_0\in a$ and $A_1\in a^\pi$. Setting $A_2:=A_1^\pi\cap a^\pi$, $a_1:=A_0\vee
A_1$ and $a_2:=A_0\vee A_2$ yields that
   \begin{displaymath}
   a\rrel a_1\rrel a_2\rrel a,\quad a^\pi\rrel a_1\rrel a_2\rrel a^\pi.
   \end{displaymath}
Hence $a_1^\varphi$ and $a_2^\varphi$ are concurrent and distinct by
(\ref{A0}). Therefore
   \begin{displaymath}
   (a_1^\varphi\cap a_2^\varphi)^\pi\cap (a_1^\varphi\vee a_2^\varphi)=:a'
   \mbox{ and } a'^\pi\neq a'
   \end{displaymath}
are the only two lines that are intersecting both $a_1^\varphi$ and
$a_2^\varphi$ orthogonally. Thus $\{a',a'^\pi\}=\{a^\varphi,a^{\pi\varphi}\}$.
Now (\ref{A3}) follows from $b^\varphi\notin\{a^\varphi,a^{\varphi\pi}\}$ and
(\ref{LR}).
\qed
\begin{prop}\label{B}
   If $a\in \Lcal$, then  either $\Sli(a)^\varphi\subset
   \Sli(a^\varphi)$  or $\Sli(a)^\varphi\subset \Sre(a^\varphi)$.
\end{prop}
\proof
Assume to the contrary that our assertion does not hold. We infer from
(\ref{LR}), (\ref{A0}), $\#\Sli(a)=\#F=\infty$ and the injectivity of
$\varphi$ that $\Sli(a^\varphi)\cap \Sre(a^\varphi) =
\{a^\varphi,a^{\varphi\pi}\} = \{a^\varphi,a^{\pi\varphi}\}$ is a proper
subset of $\Sli(a)^\varphi$. Therefore $\Sli(a)^\varphi$ cannot be a subset of
both $\Sli(a^\varphi)$ and $\Sre(a^\varphi)$. Hence there exist distinct lines
$x,y\in\Sli(a)\setminus\{a,a^\pi\}$ such that $a^\varphi\li x^\varphi$ and
$a^\varphi\re y^\varphi$. Moreover, $x^\varphi\cp y^\varphi$ by (\ref{CP}) and
(\ref{A1}).

If $x^\varphi\li y^\varphi$, then $a^\varphi\li x^\varphi\li y^\varphi$.
It follows that $y^\varphi$ is both left and right parallel to $a^\varphi$.
Hence, by (\ref{LR}) and (\ref{A2}), we obtain
$y^\varphi\in\{a^\varphi,a^{\varphi\pi}\}=\{a^\varphi,a^{\pi\varphi}\}$. This
is contradicting the injectivity of $\varphi$.

Likewise, $x^\varphi\re y^\varphi$ yields a contradiction.
\qed
\begin{prop}\label{C}
   If $\Sli(a)^\varphi\subset \Sli(a^\varphi)$ for at least one line $a\in
   \Lcal$, then
      \begin{equation}\label{C1}
      \Sli(b)^\varphi\subset \Sli(b^\varphi)
      \mbox{ for all } b\in\Lcal.
      \end{equation}
\end{prop}
\proof
At first let $a\rrel b$. Assume to the contrary that
$\Sli(b)^\varphi\not\subset \Sli(b^\varphi)$, whence Proposition \ref{B} gives
$\Sli(b)^\varphi\subset\Sre(b^\varphi)$. If $\Sli(a)$ is written down
according to (\ref{LILI-SCHNITT}), then (\ref{A0}) gives that for each
$x^\varphi\in\Sli(a)^\varphi\subset\Sli(a^\varphi)$ there exists a
$y^\varphi\in\Sli(b)^\varphi\subset\Sre(b^\varphi)$ such that $x^\varphi\rrel
y^\varphi$. We infer from (\ref{LIRE-SCHNITT}) and (\ref{RELI-SCHNITT}),
applied to $a^\varphi\rrel b^\varphi$, that
   \begin{displaymath}
   \Sli(a)^\varphi\subset\Rli{a^\varphi}{b^\varphi}
   \mbox{ and }
   \Sli(b)^\varphi\subset\Rre{b^\varphi}{a^\varphi}.
   \end{displaymath}
There exists a line $c$ such that $a\rrel c\rrel b$. Application of
(\ref{LI-RE}) yields (in an obvious shorthand notation)
$c\rrel \Rli{a}{c}$ and $c\rrel \Rli{b}{c}$. Therefore, by (\ref{A0}),
   \begin{displaymath}
   c^\varphi\rrel \Rli{a}{c}^\varphi \subset\Rli{a^\varphi}{b^\varphi}
   \mbox{ and }
   c^\varphi\rrel \Rli{b}{c}^\varphi\subset\Rre{b^\varphi}{a^\varphi}.
   \end{displaymath}
Due to the injectivity of $\varphi$,  $c^\varphi$ has to be a transversal line
of two infinite sets of lines contained in opposite reguli, respectively. This
is an absurdity.

Next assume $a\nrrel b$. Then there is a finite sequence
$a\rrel a_1\rrel\cdots\rrel a_n\rrel b$,
whence the proof from above carries over to all lines.
\qed

\begin{prop}\label{D}
   If $\Sli(a)^\varphi\subset \Sli(a^\varphi)$ for all $a\in
   \Lcal$, then
   \begin{equation}\label{D1}
      \Sre(b)^\varphi\subset \Sre(b^\varphi)
      \mbox{ for all } b\in\Lcal.
   \end{equation}
\end{prop}
\proof Given a line $b \in\Lcal$ there exists a line $a\in\Lcal$ with $b\rrel
a$. By (\ref{LI-RE}) and by the definition of the relation $\rrel$,
$\Rre{b}{a}$ and $\Rli{a}{b}$ are mutually opposite reguli containing
$\{b,b^\pi\}$ and $\{a,a^\pi\}$, respectively. There exist lines
   \begin{displaymath}
b_1\in\Rre{b}{a}\setminus\{b,b^\pi\} \mbox{ and }
a_1\in\Rli{a}{b}\setminus\{a,a^\pi\}.
   \end{displaymath}
Analogously, the distinct related lines $b^\varphi$ and $a^\varphi$ give rise
to mutually opposite reguli $\Rre{b^\varphi}{a^\varphi}$ and
$\Rli{a^\varphi}{b^\varphi}$ containing
$\{b^\varphi,b^{\pi\varphi}\} = \{b^\varphi,b^{\varphi\pi}\}$ and
$\{a^\varphi,a^{\pi\varphi}\} = \{a^\varphi,a^{\varphi\pi}\}$, respectively;
cf.\ (\ref{A2}). Now $a_1\li a$, $a_1\rrel b$, the present assumption on
$\Sli(a)^\varphi$ and (\ref{A0}) yield
   \begin{displaymath}
   a_1^\varphi\in\Rli{a^\varphi}{b^\varphi}
   \setminus\{a^\varphi,a^{\varphi\pi}\}.
   \end{displaymath}
Since $b_1\rrel \{a,a^\pi,a_1\}$, we infer from (\ref{A0}) that
$b_1^\varphi\rrel\{a^\varphi, a^{\varphi\pi},a_1^\varphi\}$. Therefore
$b_1^\varphi$ is a transversal line of the regulus
$\Rli{a^\varphi}{b^\varphi}$. So
   \begin{displaymath}
   b_1^\varphi\in
   (\Rre{b^\varphi}{a^\varphi}\setminus\{b^\varphi,b^{\varphi\pi}\})
   \subset
   (\Sre(b^\varphi)\setminus\{b^\varphi,b^{\varphi\pi}\})
   \end{displaymath}
and $b_1^\varphi\nli b^\varphi$ by (\ref{LR}). Proposition \ref{B} extends to
the $\varphi$-images of right parallel classes in an obvious way. Hence the
assertion follows.
\qed
\abstand

Propositions \ref{C} and \ref{D} hold true, mutatis mutandis, if
$\Sli(a)^\varphi\subset \Sre(a^\varphi)$. Therefore, the injection $\varphi$
is either preserving or interchanging left and right parallelism. According to
these possibilities $\varphi$ will be called {\em direct} or {\em opposite}.

In the sequel we shall confine our attention on a direct mapping $\varphi$.
The subsequent Propositions remain true when the terms ``left'' and ``right''
are interchanged.

\begin{prop}\label{E}
   If $\varphi $ is direct, then
   \begin{equation}\label{E1}
   a\nre b \mbox{ and } a\li b \Longrightarrow a^\varphi \nre b^\varphi
   \mbox{ for all } a,b\in\Lcal.
   \end{equation}
\end{prop}
\proof
We infer from $a\nre b$, $a\li b$ and (\ref{LR}) that $b\notin\{a,a^\pi\}$. By
the injectivity of $\varphi$ and (\ref{A2}), we obtain
$b^\varphi\notin\{a^\varphi,a^{\varphi\pi}\}$. Now $a^\varphi\li b^\varphi$
and (\ref{LR}) yield $a^\varphi\nre b^\varphi$.
\qed

\begin{prop}\label{F}
   If $\varphi$ is direct, then
   \begin{displaymath}
      \Abb{\varphili}{\Eli}{\Eli}{a^{\gamma\rho}}{a^{\varphi\gamma\rho}}
      \quad(a\in\Lcal)
   \end{displaymath}
   is a well-defined mapping.
\end{prop}
\proof
By (\ref{D1}), the definition of $\varphili$ is unambiguous for all points in
$\im\gamma\rho$. The restriction of $\rho$ to the Klein quadric $\Gamma$ is
surjective, since $\Gamma$ contains a plane. Thus the assertion follows.
\qed

\begin{prop}\label{H}
   Let $a\in \Lcal$. Then $\varphili|\Sli(a)^{\gamma\rho}$ is injective.
\end{prop}
\proof
Two distinct points of $\Sli(a)^{\gamma\rho}$ can be written as
$b_1^{\gamma\rho} \neq b_2^{\gamma\rho}$ with $b_1,b_2\in\Sli(a)$. Hence
$b_1\li b_2$, but $b_1\nre b_2$. By (\ref{E1}), $b_1^\varphi\nre
b_2^\varphi$. Now $b_1^{\varphi\gamma\rho} \neq b_2^{\varphi\gamma\rho}$
follows from the definition of $\re$.
\qed
\abstand

We remark that $\Sli(a)^{\gamma\rho}$ is in general a proper subset of $\Eli$:
The image of $\Sli(a)$ under the Klein mapping $\gamma$ is an oval quadric in
the 3-dimensional subspace $\Tcal:=\Eli\vee a^\gamma$ of $(\Phat,\Lhat)$; cf.\
the proof of \cite[Lemma 3]{Ha94}. The restriction of the projection $\rho$ to
$\Sli(a)^\gamma$ may be seen as a {\em gnomonic projection}, since
$\Ere\cap\Tcal$ (the centre of the projection) is an interior point of that
oval quadric. The surjectivity of this gnomonic projection is equivalent to
the fact that every right parallel class has non-empty intersection with
$\Sli(a)$. We state two sufficient conditions for
$\Sli(a)^{\gamma\rho}=\Eli$: If $F$ is a Euclidean field, then any line
through an interior point of an oval quadric is a secant; cf., for example,
\cite[vol.\ II, p.\ 54]{Br76}. Hence any gnomonic projection is surjective. If
$F$ is a Pythagorean field and if the absolute polarity $\pi$ can be described
by the standard bilinear form on $F^4$, then a line joining any point of
$\Eli$ with any point of $\Ere$ is a secant of the Klein quadric; cf.\
\cite[Remark 4]{Ha94}. Thus the gnomonic projections arising from $\rho$ are
then surjective.

\begin{prop}\label{G}
   The mapping $\abb{\varphili}{\Eli}{\Eli}$ takes $\kappali$-conjugate points
   to $\kappali$-conjugate points. Moreover, $\varphili$ is a full lineation.
\end{prop}
\proof
Let $\Fcal\subset\Lcal$ be a ruled plane. The restriction of $\gamma\rho$ to
$\Fcal$ is a collineation of $\Fcal$ (regarded as dual projective plane) onto
$\Ecal_L$. Given points $U,V\in\Eli$ there exist lines $u,v\in\Fcal$ with
$u^{\gamma\rho}=U$ and $v^{\gamma\rho}=V$. The set $\{x\mid x\in \Fcal, x\rrel
u\}$ is a pencil of lines. Its image under $\gamma\rho$ is a line. More
precisely, this is the polar line of $U$ with respect to $\kappali$; see
(\ref{RREL}). It follows from (\ref{A0}) and (\ref{RREL}) that
   \begin{equation}\label{G0}
   U,V \mbox{ $\kappali$-conjugate }
   \Longleftrightarrow
   u\rrel v
   \Longrightarrow
   a^\varphi\rrel v^\varphi
   \Longrightarrow
   U^\varphili, V^\varphili \mbox{ $\kappali$-conjugate.}
   \end{equation}

If $g\subset\Eli$ is a line, then (\ref{G0}) implies $g^\varphili\subset
g^{\kappali\varphili\kappali}$. From this observation it is immediate that
$\varphili$ is a lineation, that is, a collinearity-preserving mapping.

Finally, choose lines $a,b,c\in\Fcal$ with $a\rrel b\rrel c\rrel a$.
By (\ref{RREL}), $\{a^{\gamma\rho},b^{\gamma\rho},c^{\gamma\rho}\}
\subset\Eli$ is a self-polar triangle with respect to $\kappali$. We form the
reguli
   \begin{equation}\label{G1}
   \Rli{c}{a}\rrel a,\; \Rli{a}{b}\rrel b,\;\Rli{b}{c}\rrel c.
   \end{equation}
If $\Rcal$ denotes any of these three reguli, then the restriction
$\gamma\rho|\Rcal$ is not injective, but the fibre of $x\in\Rcal$ is given by
$\{x,x^\pi\}$. There is an infinite number of such unordered pairs, whence the
image $\Rcal^{\gamma\rho}$ is an infinite subset on one side of the triangle
$\{a^{\gamma\rho},b^{\gamma\rho},c^{\gamma\rho}\}$; e.g.,
$\Rli{c}{a}^{\gamma\rho}\subset(b^{\gamma\rho}\vee c^{\gamma\rho})$. By
(\ref{G0}) and Proposition \ref{H}, $\im\varphili$ contains the
$\kappali$-self-polar triangle
$\{a^{\varphi\gamma\rho},b^{\varphi\gamma\rho},c^{\varphi\gamma\rho}\}$ as
well as infinitely many points on each side of this triangle. We read off from
\cite[p.\ 4]{CV80} that $\im\varphili$ is either a projective subplane or a
near-pencil (degenerate subplane) of $\Eli$. Obviously, the second possibility
cannot occur, whence there exists a quadrangle in $\im\varphili$, i.e.,
$\varphi_L$ is full.
\qed
\abstand

The next result is immediate from Proposition \ref{H}, whenever we have a
surjective gnomonic projection $\rho|\Sli(a)^\gamma$ for at least one line
$a\in \Lcal$.

\begin{prop}\label{I}
The lineation $\abb{\varphili}{\Eli}{\Eli}$ injective.
\end{prop}
\proof
Choose intersecting lines $a,t\in\Lcal$ with $a\nrel t$. Set $\Tcal:=
\spn\Sli(a)^\gamma$. The regulus $\Rli{a}{t}$ does not contain $a^\pi$ so that
the plane
{\unitlength1cm
   \begin{figure}[h]
      \begin{center}
         \begin {picture}(10.9,8.0)
         \put(0.3,0.1){\includegraphics[width=10.9\unitlength]{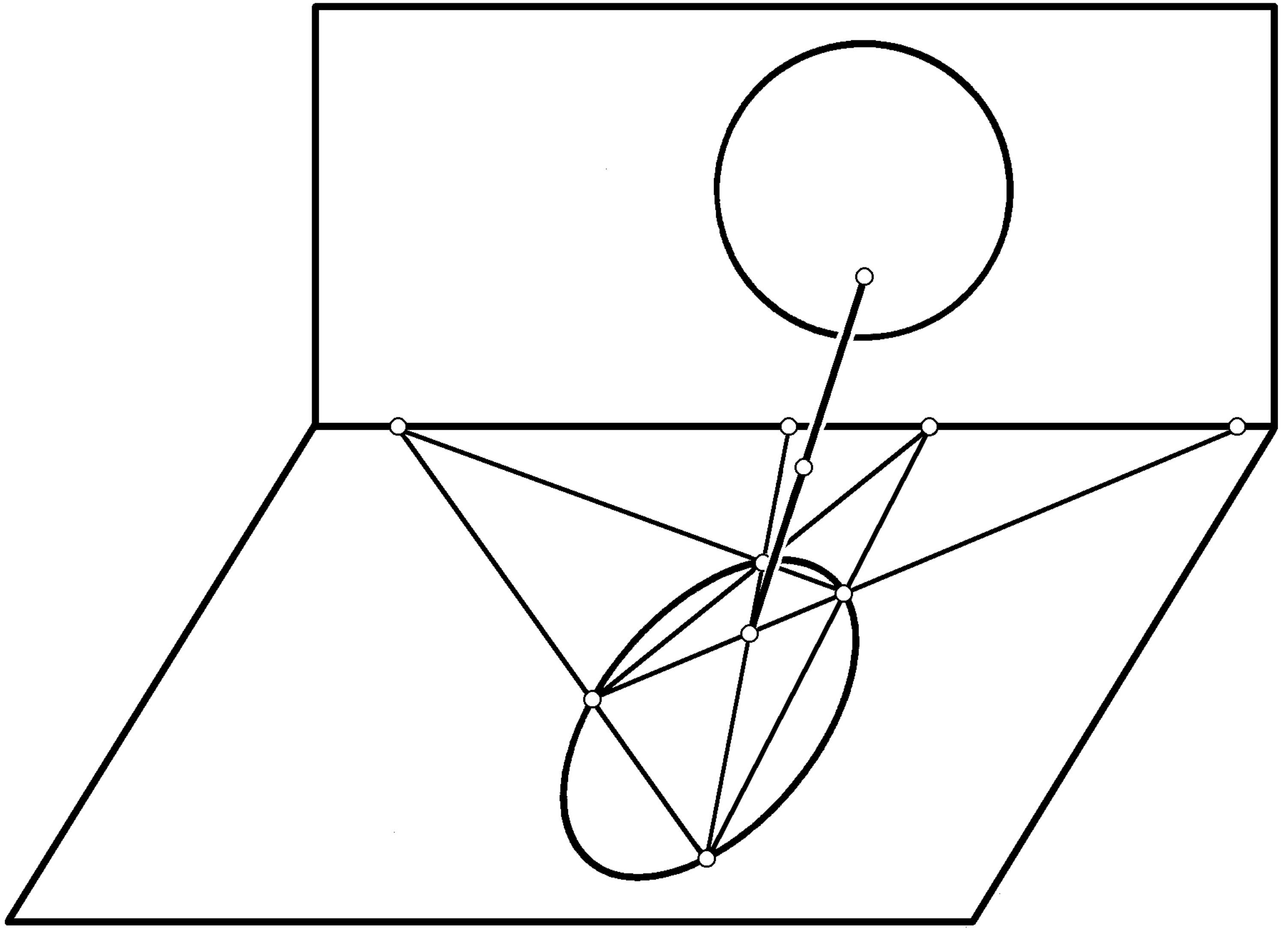}}
         \put(6.3,3.4){$A$}
         \put(6.3,0.25){$B$}
         \put(6.8,4.5){$C$}
         \put(7.5,1.7){$k$}
         \put(6.7,2.2){$U$}
         \put(5.0,4.5){$u$}
         \put(4.75,1.8){$X$}
         \put(7.7,2.6){$\qu{X}$}
         \put(8.0,4.5){$X_A$}
         \put(3.5,4.5){$X_B$}
         \put(9.0,3.2){$C^\kappali$}
         \put(0.9,0.5){$\Eli$}
         \put(3.2,7.5){$\Wcal$}
         \put(9.0,6.3){$\Rli{a}{t}^\gamma$}
         \put(7.45,4.7){$m$}
         \put(7.15,3.7){$\Ere\cap\Tcal$}

         \end{picture}
      \end{center}
   \end{figure}
}%
$\Wcal:=\spn\Rli{a}{t}^\gamma\subset \Tcal$ does not contain the point
$\Ere\cap\Tcal$; cf.\ the figure. Therefore the image of the conic
$\Rli{a}{t}^\gamma$ under the gnomonic projection $\rho|\Sli(a)^\gamma$, is
also a conic, say
   \begin{equation}\label{I00}
   k:=\Rli{a}{t}^{\gamma\rho}.
   \end{equation}
Put $u:=\Wcal\cap\Eli$ and write $m$ for the polar line of $u$ with respect to
the polarity of $\Sli(a)^\gamma$. Then $m\cap\Wcal$ is the pole of $u$ with
respect to $\Rli{a}{t}^\gamma$. As$\Ere\cap\Tcal$ is on $m$, the point
$U:=m\cap\Eli=(m\cap\Wcal)^\rho$ is the pole of $u$ with respect to $k$ and
$\kappali$. The elliptic polarity
$\kappali$ and the polarity of the conic $\Rli{a}{t}^\gamma$ induce the same
elliptic involutory projectivity on $u$, since both are arising from the
polarity $\kappa$ of the Klein quadric. The gnomonic projection fixes $u$
pointwise, whence this involutory projectivity on $u$ is also induced by the
polarity of the conic%
   \footnote{$k$ is a circle of the elliptic plane $\Eli.$}
$k$. Hence $U$ is an interior point of $k$. By Proposition \ref{H} and
(\ref{I00}),
   \begin{equation}\label{I0}
\varphili|k \mbox{ is injective.}
   \end{equation}
Since $k$ is infinite, we may find points $A,B\in k$, collinear with $U$, such
that $A^\varphili$, $B^\varphili$ and $U^\varphili$ are mutually distinct. By
Proposition \ref{G},
   \begin{equation}\label{I1}
   U^\varphili\in A^\varphili\vee B^\varphili.
   \end{equation}

Let $X\in k\setminus\{A,B\}$. As $U$ is interior point of $k$, the line $X\vee
U$ meets $k$ at one more point $\qu{X}\neq X$. Then $\{A,B,X,\qu{X}\} \subset
k$ is a quadrangle, whence its diagonal points form a self-polar triangle with
respect to $k$, say $\{U,X_A,X_B\}$. Hence $u=X_A\vee X_B$,
and $\{U,X_A,X_B\}$ is also a $\kappali$-self-polar triangle. We infer from
(\ref{I0}) that $A^\varphili,B^\varphili,X^\varphili,\qu{X}^\varphili$ are
mutually distinct. Proposition \ref{G} tells us that
$\{U^\varphili,X_A^\varphili,X_B^\varphili\}$ is a $\kappali$-self-polar
triangle. We claim that
   \begin{equation}\label{I2}
   X^\varphili\notin A^\varphili\vee B^\varphili;
   \end{equation}
otherwise we would have
   \begin{displaymath}
   X_A^\varphili\in A^\varphili\vee X^\varphili = A^\varphili\vee B^\varphili
\mbox{ and }
   X_B^\varphili\in B^\varphili\vee X^\varphili = A^\varphili\vee B^\varphili,
   \end{displaymath}
whence $X_A^\varphili$, $X_B^\varphili$ and $U^\varphili$ would be collinear,
an absurdity.

Since $X\in k\setminus\{A,B\}$ has been chosen arbitrarily, we can deduce from
(\ref{I1}) and (\ref{I2}) that
   \begin{equation}\label{I3}
   P^\varphili\neq U^\varphili \mbox{ for all } P\in k.
   \end{equation}
It is obvious now that $A^\varphili,B^\varphili,X^\varphili,\qu{X}^\varphili$
is a quadrangle with $U^\varphili$ being one of its diagonal points. Thus we
have arrived at the well-known construction of the fourth harmonic point for
$A^\varphili,B^\varphili,U^\varphili$. Setting
   \begin{displaymath}
      u':=U^{\varphili\kappali} = X_A^\varphili\vee X_B^\varphili
   \end{displaymath}
yields that
   \begin{equation}
   u'\neq A^\varphili\vee B^\varphili.
   \end{equation}
With $C:=u\cap (A\vee B)$, we obtain $C^\varphili=u'\cap (A^\varphili \vee
B^\varphili)$ and
   \begin{equation}\label{I4}
   C^\varphili \neq X_A^\varphili.
   \end{equation}

As $X$ varies in $k\setminus\{A,B\}$, the point $X_A$ is running in
$u\setminus(\{C\}\cup (C^\kappali\cap u))$, that is, we are reaching all
points of $u$ but two. By (\ref{I4}), $C^\varphili\in u'$ has only a finite
number of pre-images on the line $u$. Now \cite[Satz 3.2]{An70} establishes
that $\varphili$ is injective.
\qed
\begin{prop}\label{M}
   If $\abb{\varphi}{\Lcal}{\Lcal}$ is direct, then any two points of $\Eli$
   that are not $\kappali$-conjugate remain non-conjugate under $\varphili$.
\end{prop}
\proof
The injective and full lineation $\varphili$ is preserving non-collinearity of
points \cite[p.\ 4]{CV80}. If $U,V\in \Eli$ are not $\kappali$-conjugate,
then $V\notin U^\kappali$, whence $V^\varphili\notin
\spn(U^{\kappali\varphili})= U^{\varphili\kappali}$. Thus $U^\varphili,
V^\varphili$ are not $\kappali$-conjugate.
\qed
\begin{prop}\label{J}
If $\abb{\varphi}{\Lcal}{\Lcal}$ is direct, then
      \begin{equation}\label{J1}
      a\nrel b \Longrightarrow a^\varphi \nrel b^\varphi \mbox{ for all }
      a,b\in\Lcal.
      \end{equation}
\end{prop}
\proof
We infer from the left hand side of (\ref{J1}) and from (\ref{RREL}) that (up
to interchanging the terms ``left'' and ``right'')
$a^{\gamma\rho},b^{\gamma\rho}$ are not $\kappali$-conjugate. Under
$\varphili$ this property remains unchanged by Proposition \ref{M}, whence the
assertion follows from (\ref{RREL}) and the definition of $\varphili$ in
Proposition \ref{F}.
\qed
\abstand

We recall a concept introduced in \cite{Ha94}: Two collineations
$\abb{\zeta}{\Eli}{\Eli}$ and $\abb{\eta}{\Ere}{\Ere}$ are called {\em
admissible} if the following conditions hold true:

\abstand
   \begin{tabular}{ll}
   Ad1. &
$\zeta \mbox{ and }\eta \mbox{ are commuting with }\kappali \mbox{
   and }\kappare,\mbox{ respectively}.$\\
   Ad2. &
   $(X\vee Y)\cap\Gamma\neq\emptyset \Longrightarrow (X^\zeta\vee
   Y^\eta)\cap\Gamma\neq\emptyset$ for all $X\in \Eli,\, Y\in\Ere$.
   \end{tabular}
\abstand

When writing \cite{Ha94}, the author considered the next Proposition to be
self-evident. It seems, however, that it deserves a formal proof.
   \begin{prop}\label{K}
   If collineations $\abb{\zeta}{\Eli}{\Eli}$ and $\abb{\eta}{\Ere}{\Ere}$ are
   admissible, then their inverse mappings are also admissible.
   \end{prop}
\proof
Condition Ad1 is evidently true for $\zeta\inv,\eta\inv$.
Given $X\in \Eli$ and $Y\in \Ere$ such that $(X\vee Y)\cap\Gamma=\emptyset$,
then, by the surjectivity of $\gamma\lambda$, there is a line $a\in \Lcal$
with $Y=a^{\gamma\lambda}$. By Ad2, there exists a line
$a'\in \Lcal$ with $a'^\gamma\in (a^{\gamma\rho\zeta}\vee Y^\eta)\cap\Gamma$.
It is straightforward to verify that $\zeta$ extends to a collineation
   \begin{displaymath}
   \abb{\delta}{Y\vee\Eli}{Y^\eta\vee\Eli}
   \end{displaymath}
with $Y\mapsto Y^\eta$ and $a^\gamma\mapsto a'^\gamma$. Under $\delta$ the
elliptic quadric $\Sli(a)^\gamma$ goes over to an elliptic quadric within
$Y^\eta\vee\Eli$. This quadric coincides with the quadric
$(Y^\eta\vee\Eli)\cap\Gamma$, since $a'^\gamma$ is a common point, $\Eli$ is
the common polar plane of $Y^\eta$ and $\kappali=\zeta\inv\kappali\zeta$ is
the induced polarity in $\Eli$ for both quadrics%
   \footnote{This is a projective generalization of the well-known fact that a
   sphere in Euclidean space is uniquely determined by one point and its
   mid-point.}%
; cf., e.g.,
\cite[vol.\ I, p.\ 191]{Br76}. Thus $(X^\zeta\vee Y^\eta)\cap\Gamma = ((X\vee
Y)\cap\Gamma)^\delta=\emptyset$.
\qed
\abstand

{\it Proof of Theorem \ref{TH1}. }
If $\varphi$ is direct, then the result is immediate from (\ref{J1}).
Otherwise
choose a point $Q\in\Pcal$. The harmonic homology with centre $Q$ and axis
$Q^\pi$ is an elliptic reflection and yields an opposite Pl\"ucker
transformation $\chi$. By (\ref{J1}), $\chi\varphi$ is a direct Pl\"ucker
transformation so that $\varphi$ too is a Pl\"ucker transformation.
\qed
\abstand

{\it Proof of Theorem \ref{TH2}. }
If $\varphi$ is direct, then $\Eli^\varphili$ is a subplane of $\Eli$
isomorphic to $\Eli$. Hence the underlying field of $\Eli^\varphili$ is a
subfield of $F$ isomorphic to $F$ \cite[p.\ 266]{St72}. This implies, by our
assumption on $F$, that $\varphili$ is surjective. We infer from Propositions
\ref{G} and \ref{I} that $\varphili$ is a collineation of $\Eli$ commuting
with $\kappali$. Similarly, $\varphire$ is a collineation of $\Ere$ commuting
with $\kappare$. By their definition, $\varphili$ and $\varphire$ are
admissible.

Given a line $a'\in \Lcal$ we may apply $\varphiliinv$ and $\varphireinv$ to
the points $a'^{\gamma\rho}$ and $a'^{\gamma\lambda}$, respectively. Thus we
obtain points $X\in\Eli$ and $Y\in\Ere$, say. By Proposition \ref{K}, there
exists a line $a\in \Lcal$ such that $a^\gamma$ is on the line $X\vee
Y$. Therefore $\{a,a^\pi\}^\varphi=\{a',a'^\pi\}$ so that $\varphi$ is
surjective. The assertion follows now from Theorem \ref{TH1}.

If $\varphi$ is opposite, then let $\chi$ be an elliptic reflection. Hence
both $\chi\varphi$ and $\varphi$ are Pl\"ucker transformations.
\qed


\abstand

\noindent
Hans Havlicek\\
Abteilung f\"ur Lineare Algebra und Geometrie\\
Technische Universit\"at\\
Wiedner Hauptstra{\ss}e 8--10\\
A-1040 Wien, Austria\\
EMAIL: {\tt havlicek@geometrie.tuwien.ac.at}
\end{document}